\documentclass[12pt]{article}
\newtheorem{theorem}{Theorem}[section]
\newtheorem{lemma}{Lemma}[section]
\newtheorem{definition}{Definition}[section]

\newtheorem{corollary}{Corollary}[section]
\newtheorem{proposition}{Proposition}[section]
\usepackage{setspace}
\usepackage{mathrsfs}
\usepackage{amsmath}
\usepackage{amssymb}
\usepackage{amsfonts}
\usepackage{amsmath,amsfonts,amssymb}
\usepackage{graphicx}
\usepackage{subfigure}
\usepackage{caption2}
\usepackage[rightcaption]{sidecap}
\usepackage{fancyhdr}
\newcommand{\R}{\mathbb R}
\newcommand{\B}{\mathbb B}

\newtheorem{remark}{Remark}

\begin{document}
\title{Properties of the extremal solution for a fourth-order elliptic problem
}
\date{}
\author{Baishun Lai$^{1,2}$, Zhuoran Du\\
\small {\it 1 Institute of Contemporary Mathematics, Henan University}\\
\small {\it 2  School of Mathematics and Information Science,Henan University}\\
\small { Kaifeng 475004, P.R.China}\\
\small {\it College of Mathematics and Econometrics, Hunan University}\\
\small { Changsha 410082, P.R. China}\\
}
\maketitle
\begin{center}
\begin{minipage}{130mm}
{\small {\bf Abstract}\ \ \  Let $\lambda^{*}>0$ denote the largest possible value of $\lambda$ such that
$$
\left\{
\begin{array}{lllllll}
\Delta^{2}u=\frac{\lambda}{(1-u)^{p}} &  \mbox{in}\ \ \B, \\
0<u\leq 1 &  \mbox{in}\ \ \B, \\
u=\frac{\partial u}{\partial n} =0 &  \mbox{on}\  \  \partial \B\\
\end{array}
\right.
$$
has a solution, where $\B$ is the unit ball in $\R^{n}$ centered at
the origin, $p>1$ and $n$ is the exterior unit normal vector. We
show that for $\lambda=\lambda^{*}$ this problem possesses a unique
weak solution $u^{*}$, called the extremal solution. We prove that
$u^{*}$ is singular when $n\geq 13$ for $p$ large enough and
 actually solve part of the open problem which
\cite{D} left.\vskip 0.3in
{\bf Mathematics Subject Classification (2000)} \  \  35B45 $\cdot$ 35J40.
}
\end{minipage}
\end{center}
\vskip 0.2in

\setcounter{equation}{0}
 \setcounter{section}{0}
\section{Introduction and result}

 The main purpose of this paper is to investigate
regularity of the extremal solution for a class of fourth-order
problem
$$
\left\{
\begin{array}{lllllll}
\Delta^{2}u=\frac{\lambda}{(1-u)^{p}} &  \mbox{in}\ \ \B, \\
0<u\leq 1 &  \mbox{in}\ \ \B, \\
u=\frac{\partial u}{\partial n} =0 &  \mbox{on}\  \  \partial \B.\\
\end{array}
\right.
\eqno(1.1)_{\lambda}
$$
 Here $\B$ denotes the unit ball in $\R^{n}$ ($n\geq 2$) centered at the origin, $\lambda>0, p>1$ and $\frac{\partial }{\partial n}$
 the differentiation with the respect to the exterior unit normal, i.e., in radial direction. We consider only radial solutions,
since all positive smooth solutions of $(1.1)_{\lambda}$ are radial, see Berchio et al. \cite{Be}.\vskip 0.1in

The motivation for studying $(1.1)_{\lambda}$ stems from a model for the steady sates of a simple micro electromechanical
system (MEMS) which has the general form (see for example  \cite{Li}, \cite{Pe})
$$
\left\{
\begin{array}{lllllll}
\alpha \Delta^{2}u=(\beta\int_{\Omega}|\nabla u|^{2} dx+\gamma)\Delta u+
\frac{\lambda f(x)}{(1-u)^{2}(1+\chi\int_{\Omega}\frac{dx}{(1-u)^{2}})}& \mbox{in}\ \ \Omega,\\
0<u<1&  \mbox{in}\ \ \Omega,\\
u=\alpha\frac{\partial u}{\partial n}=0 &  \mbox{on }\ \partial\Omega,
\end{array}
\right.
\eqno(1.2)
$$
where $\alpha, \beta,\gamma,\chi\geq0,$ are fixed, $f\geq 0$ represents the permittivity profile,
$\Omega$ is a bounded domain in $\R^{n}$ and $\lambda>0$ is a constant which is increasing with respect to the applied voltage.\vskip 0.1in

Recently, Equation (1.2) posed in $\Omega=\B$ with $\beta=\gamma=\chi=0, \alpha=1$ and $f(x)\equiv1$, which is reduced to
$$
\left\{
\begin{array}{lllllll}
\Delta^{2}u=\frac{\lambda}{(1-u)^{2}} & \mbox{in}\ \ \B,\\
0<u<1&  \mbox{in}\ \ \B, \\
u=\frac{\partial u}{\partial n}=0 &  \mbox{on } \ \partial \B,
\end{array}
\right.
\eqno(1.3)
$$
has been studied extensively in \cite{Co1}. For convenience, we now give the following notion of solution.\vskip 0.1in

\begin{definition}\label{D1.1}
If $u_{\lambda}$ is a solution of $(1.1)_{\lambda}$ such that for
any other solution $v_{\lambda}$ of $(1.1)_{\lambda}$ one has
$$
u_{\lambda}\leq v_{\lambda}, \ \ \ \mbox{ a.e.} \ x\ \in \B,
$$
we say that $u_{\lambda}$ is a minimal solution of $(1.1)_{\lambda}$.
\end{definition}

It is shown that there exists a critical value $\lambda^{*}>0$ (pull-in voltage) such that if $\lambda\in (0,\lambda^{*})$
the problem (1.3) has a smooth minimal solution , while for $\lambda>\lambda^{*}$
(1.3) has no solution even in a weak sense. Moreover, the branch $\lambda\to u_{\lambda}(x)$ is increasing for each $x\in \B$,
and therefore the function $u^{*}(x):=\lim_{\lambda\to\lambda^{*}}u_{\lambda}(x)$ can be considered as
a generalized solution that corresponds to the pull-in voltage $\lambda^{*}$. Now the issue of the regularity of this extremal
solution-which, by elliptic regularity theory, is equivalent to whether $\sup_{\B}u^{*}<1$- is an important question for many reasons.
 For example, one of the reason is that it decides whether the set of solutions stops there, or whether a new branch of solutions emanates
 from a bifurcation state $(u^{*}, \lambda^{*})$ (see Figures 1,2). This issue turned out to depend closely on the dimension. Indeed by the
 key uniform estimate of $\|(1 - u)^{-3}\|_{L^1}$, Guo and Wei \cite{Gu4} obtained the regularity of the extremal solution
  for small dimensions and
 they proved that for dimension $n=2$ or $n=3, u^{*}$ is smooth. But from their result, the regularity of extremal solution of $(1.3)$
 is unknown for $n\geq 4$.  Recently, using certain improved Hardy-Rellich inequalities, Cowan-Esposito-Ghoussoub-Moradfam \cite{Co1} improved the above result
 and they obtained that $u^{*}$ is regular in dimensions $1\leq n\leq 8$, while it is singular for $n\geq 9$, i.e., the critical dimension is 9.
So the issue of the regularity of the extremal solution of
$(1.1)_{\lambda}$ for power $p=2$ is completely solved, but the
critical dimension for generally  power is {\bf unknown}.
\vskip0.1in
\begin{figure}
\includegraphics[width=5cm,height=4cm]{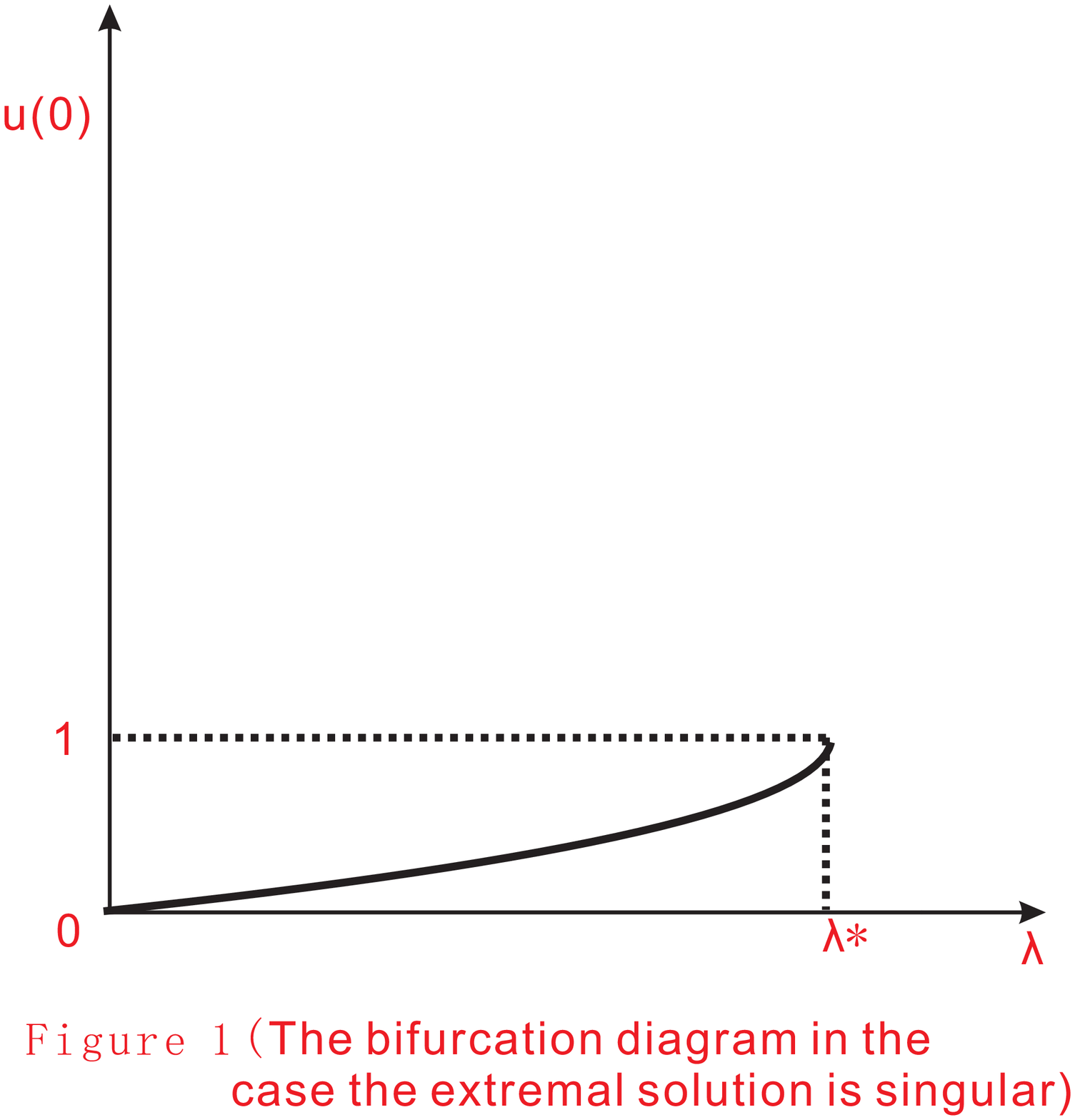}
\hspace{2.5cm}
\includegraphics[width=5cm,height=4cm]{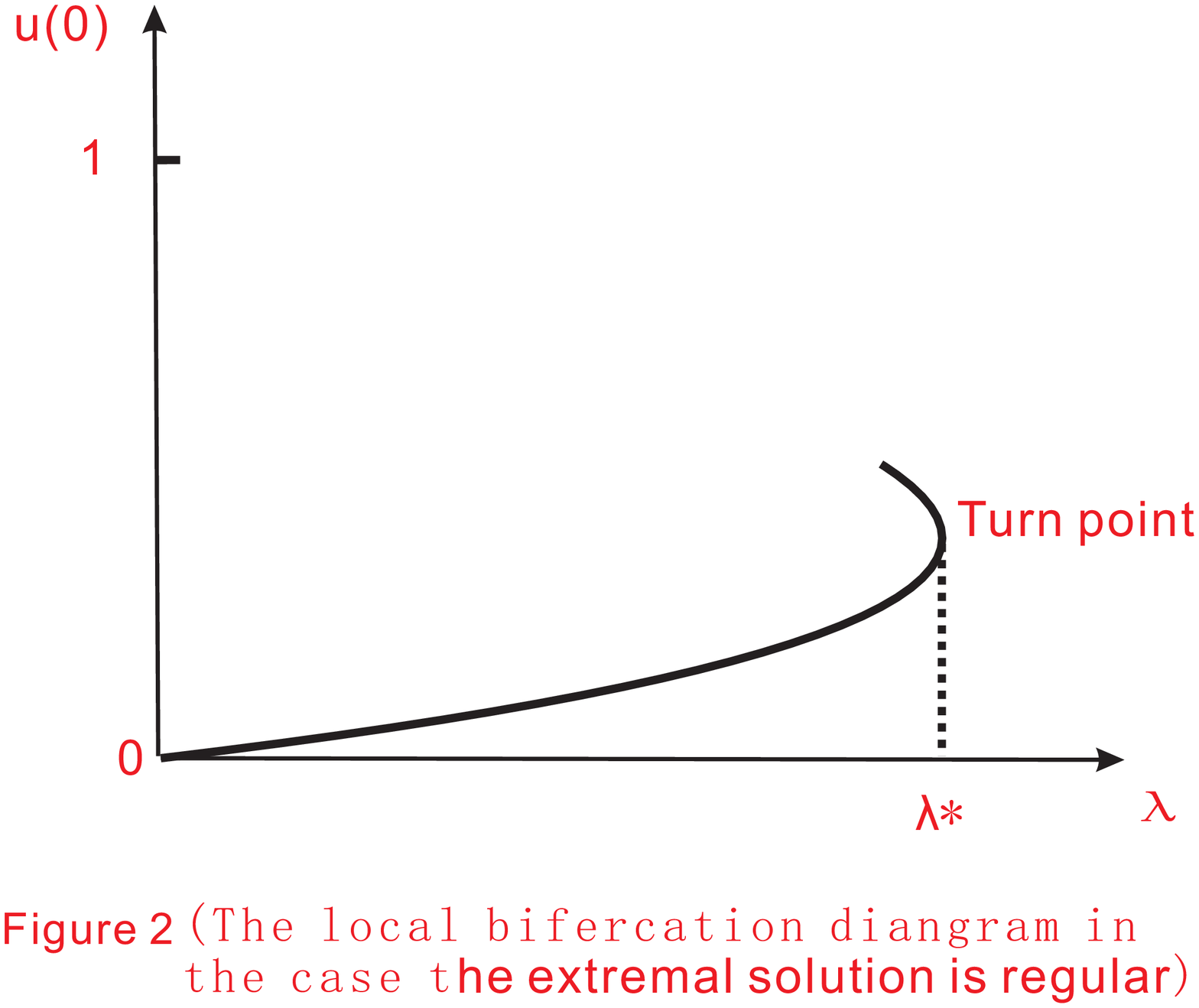}
\end{figure}

Recently, the multiplicity phenomenon for radial solutions of
$(1.1)_{\lambda}$ and the regularity of the extremal solution of
$(1.1)_{\lambda}$ for a large range of  powers have been studied
extensively by Juan D\`{a}vila etal \cite{D}.  For convenience, we
now define:
$$
p_{c}=\frac{n+2-\sqrt{4+n^{2}-4\sqrt{n^{2}+H_{n}}}}{n-6-\sqrt{4+n^{2}-4\sqrt{n^{2}+H_{n}}}} \ \ \quad \mbox{for}\ \  n\geq 3;
$$
$$
p_{c}^{+}=\frac{n+2+\sqrt{4+n^{2}-4\sqrt{n^{2}+H_{n}}}}{n-6+\sqrt{4+n^{2}-4\sqrt{n^{2}+H_{n}}}} \ \  \quad  \mbox{for}\ \  n\geq 3,\ \ n\neq 4
$$
with $H_{n}=(n(n-4)/4)^{2}$ and the numbers $p_{c}$ and $p_{c}^{+}$ are  such that when $-p=p_{c}$ or $-p=p_{c}^{+}$ then
$$
(\frac{4}{-p-1}+4)(\frac{4}{-p-1}+2)(n-2-\frac{4}{-p-1})(n-4-\frac{4}{-p-1})=H_{n}.
$$
\vskip0.1in

To explain our motivations, we now recall some corresponding results
from $\cite{D}$\vskip0.1in

\noindent \textit{{\bf  Theorem A} \ Assume
$$
n=3 \ \ \mbox{and}\ \  p_{c}^{+}<-p<p_{c},\ \ \mbox{or}\ \ 4\leq n\leq 12\ \ \mbox{and}\ \ -\infty<-p<p_{c}.\eqno(1.5)
$$
Then there exist a unique $\lambda_{s}$ such that  $(1.1)_{\lambda}$
with $\lambda=\lambda_{s}$  has infinitely many radial smooth
solutions. For $\lambda\neq\lambda_{s}$ there are finitely many
radial smooth solutions and their number goes to infinity as
$\lambda\to \lambda_{s}$. Moreover, $\lambda_{s}<\lambda^{*}$ and
$u^{*}$ is regular.}\vskip0.1in
From this Theorem,  we know that the
extremal solution of $(1.1)_{\lambda}$ is regular for a certain
range of $p$ and $n$.
 At the same time, they left a open problem:
if
$$
\left\{
\begin{array}{lllllll}
n=3\  \mbox{and}\  -p\in (-3,p_{c}^{+}]\cup [p_{c},-1)\  \mbox{or}\\
5\leq n\leq 12 \ \ \mbox{and}\ \ p_{c}\leq -p<-1, \ \mbox{or}\\
n\geq 13 \ \ \mbox{and}\ -p<-1,
\end{array}
\right.
$$
is $u^{*}$ singular?

In this paper, by constructing a semi-stable singular $H_{0}^{2}(\B)-$ weak sub-solution of $(1.1)_{\lambda}$, we prove that, if $p$ is large enough,
 the extremal solution is singular for dimensions $n\geq 13$ and complete part of the above open problem. Our result is stated as follows:
\vskip 0.2in

\begin{theorem}\label{T1.1}
(i) For any $p>1$, the  unique extremal solution of
$(1.1)_{\lambda^{*}}$ is regular for dimensions $n\leq 4$;\vskip
0.1in \noindent
 (ii)  There exists $p_0 > 1$ large enough such that
for $p \geq p_0$, the unique extremal solution of
$(1.1)_{\lambda^{*}}$ is singular for dimensions $n\geq13$.\vskip
0.2in
\end{theorem}

From the technical point of view, one of the basic tools in the
analysis of nonlinear second order elliptic problems in bounded and
unbounded domains
 of $\R^{n} (n\geq 2)$ is the maximum principle. However, for high order problems, such principle dose not normally hold for general domains
 (at least for the clamped boundary conditions $u=\frac{\partial u}{\partial n}=0$ on $\partial \Omega$), which causes several technical difficulties.
 One of reasons to the study $(1.1)_{\lambda}$ in a ball is that a maximum principle holds in this situation, see \cite{Ar1}, \cite{Bo}.
 The second obstacle is the well-known difficulty of extracting energy estimates for solutions of fourth order problems from their stability properties.
  Besides, for the corresponding second order problem, the starting point was an explicit singular solution for a suitable eigenvalue
parameter $\lambda$ which turned out to  play a fundamental role for the shape of the corresponding bifurcation diagram, see \cite{Br}.
When turning to the biharmonic problem $(1.1)_{\lambda}$ the second boundary condition $\frac{\partial u}{\partial n}=0$
prevents to find an explicit singular solution.
This means that the method used to analyze the regularity of the extremal solution for second order problem could not carry to the
 corresponding problem for $(1.1)_{\lambda}$. In this paper, we, in order to overcome the third obstacle, use improved and non
standard Hardy-Rellich inequalities recently established by Ghoussoub-Moradifam in \cite{Gh} to construct a semi-stable singular
 $H^{2}(\B)-$ weak sub-solution of $(1.1)_{\lambda}$.
\vskip 0.1in

This paper is organized as follows. In the next section, some preliminaries are reviewed.  In Section 3, we give the uniform estimate of
$\|(1-u)^{-(p+1)}\|_{L^{1}}$ according to the stability of the minimal solutions.
 We study the regularity of the extremal solution of $(1.1)_{\lambda}$ and the Theorem 1 (i) is established in Section 4. Finally,
  we will show that the extremal solution $u^{*}$ in dimensions $n\geq 13$ is singular by constructing a semi-stable singular $H^{2}(\B)-$
  weak sub-solution of $(1.1)_{\lambda}$.

\setcounter{equation}{0}
 \setcounter{section}{1}
\section{Preliminaries}
First we give some comparison principles which will be used throughout the paper\vskip0.1in

\begin{lemma}\label{L2.1}
(Boggio's principle, \cite{Bo}) If $u\in C^{4}(\bar{\B}_{R})$
satisfies
$$
\left\{
\begin{array}{lllllll}
\Delta^{2}u\geq 0 & \mbox{in}\ \  \B_{R},\\
u=\frac{\partial u}{\partial n}=0 &  \mbox{on}\ \  \partial \B_{R},
\end{array}
\right.
$$
then $u\geq 0$ in $\B_{R}$. \vskip0.1in
\end{lemma}

\begin{lemma}\label{L2.2}
Let $u\in L^{1}(\B_{R})$ and suppose that
$$
\int_{\B_{R}}u\Delta^{2}\varphi\geq0
$$
for all $\varphi\in C^{4}(\bar{\B}_{R})$ such that $\varphi\geq0$ in
$\B_{R}$, $\varphi|_{\partial \B_{R}}=\frac{\partial
\varphi}{\partial n}|_{\partial \B_{R}}=0$. Then $u\geq 0$ in
$\B_{R}$. Moreover $u\equiv 0$ or $u>0$ a.e., in $\B_{R}$.
\end{lemma}
For a proof see Lemma 17 in \cite{Ar1}. \vskip0.1in

\begin{lemma}\label{L2.3}
 If $u\in H^{2}(\B_{R})$ is radial, $\Delta^{2}u\geq 0$ in $\B_{R}$ in the weak sense, that is
$$
\int_{\B_{R}}\Delta u\Delta\varphi\geq 0 \ \ \forall \varphi \in C_{0}^{\infty}(\B_{R}), \ \varphi\geq0
$$
and $u|_{\partial \B_{R}}\geq0, \frac{\partial u}{\partial n}|_{\partial \B_{R}}\leq 0$ then $u\geq 0$ in $\B_{R}$.\vskip0.1in
\end{lemma}
{\bf Proof.}  We only deal with the case $R=1$ for simplicity. Solve
$$
\left\{
\begin{array}{lllllll}
\Delta^{2}u_{1}=\Delta^{2}u & \mbox{in} \ \  \B\\
u_{1}=\frac{\partial u_{1}}{\partial n}=0 &  \mbox{on}\ \ \partial \B
\end{array}
\right.
$$
in the sense $u_{1}\in H_{0}^{2}(\B)$ and $\int_{\B}\Delta u_{1}\Delta\varphi=\int_{\B}\Delta u\Delta\varphi$
for all $\varphi \in C_{0}^{\infty}(\B)$. Then $u_{1}\geq 0$ in $\B$ by lemma 2.2.

Let $u_{2}=u-u_{1}$ so that $\Delta^{2}u_{2}=0$ in $\B$. Define $f=\Delta u_{2}$. Then $\Delta f=0$ in $\B$ and since $f$
is radial we find that $f$ is a constant. It follows that $u_{2}=ar^{2}+b$. Using the boundary conditions we deduce $a+b\geq 0$
 and $a\leq0$, which imply $u_{2}\geq0$.\vskip 0.1in

As in \cite{Co1}, we are now led here to examine problem $(1.1)_{\lambda}$ with non-homogeneous boundary conditions such as

$$
\left\{
\begin{array}{lllllll}
\Delta^{2}u=\frac{\lambda}{(1-u)^{p}}&\ \ \ \mbox{in}\ \B,\\
\alpha<u\leq 1& \ \ \  \mbox{in}\ \B,\\
u=\alpha, \frac{\partial u}{\partial n}=\gamma & \ \ \  \mbox{on}\ \partial \B,
\end{array}
\right.
\eqno(2.1)_{\lambda,\alpha,\gamma}
$$
where $\alpha,\gamma$ are given.\vskip0.1in

Let $\Phi$ denote the unique solution of
$$
\left\{
\begin{array}{lllllll}
 \Delta^{2}\Phi=0&\ \ \ \mbox{in}\  \B,\\
\Phi=\alpha, \frac{\partial \Phi}{\partial n}=\gamma&\ \ \ \mbox{on}\ \partial \B.
\end{array}
\right.
\eqno(2.2)
$$
We will say that the pair $(\alpha,\gamma)$ is admissible if $\gamma\leq0$, and $\alpha-\frac{\gamma}{2}<1$.
 We now introduce a notion of weak solution. \vskip0.1in

\begin{definition}\label{D2.1}
We say that $u$ is a weak solution of
$(2.1)_{\lambda,\alpha,\gamma}$, if $\alpha\leq u\leq1$ a.e. in
$\Omega$, $\frac{1}{(1-u)^{p}} \in L^{1}(\Omega)$ and if
$$
\int_{\B}(u-\Phi)\Delta^{2}\varphi=\lambda \int_{\B}\frac{\varphi}{(1-u)^{p}}\ \ \
\forall \varphi \in C^{4}(\bar{\B})\cap H_{0}^{2}(\B),
$$
where $\Phi$ is given in (2.2). We say $u$ is a weak super-solution (resp. weak sub-solution) of  $(2.1)_{\lambda,\alpha,\gamma}$,
if the equality is replaced with $\geq$ (resp.$\leq$) for $\varphi\geq 0$.
\end{definition}

\begin{definition}\label{D2.2}
 We say a weak solution of $(2.1)_{\lambda,\alpha,\gamma}$ is regular (resp. singular)
 if $ \|u\|_{\infty}<1$ (resp. $\|u\|=1$) and stable (resp. semi-stable) if
$$
\mu_{1}(u)=\inf\{\int_{\B}(\Delta \varphi)^{2}-p\lambda\int_{\B}\frac{\varphi^{2}}{(1-u)^{p+1}}:
\phi\in H_{0}^{2}(\B), \|\phi\|_{L^{2}}=1\}
$$
is positive (resp. non-negative).
\end{definition}

We now define
$$
\lambda^{*}(\alpha,\gamma):=\sup\{\lambda>0: (2.1)_{\lambda,\alpha,\gamma}\ \mbox{has a classical soltion}\}
$$
and

$$
\lambda_{*}(\alpha,\gamma):=\sup\{\lambda>0: (2.1)_{\lambda,\alpha,\gamma}\ \mbox{has a weak soltion}\}.
$$
Observe that by Implicit Function Theorem, we can classically solve $(2.1)_{\lambda,\alpha,\gamma}$ for small $\lambda's$.
Therefore, $\lambda^{*}(\alpha,\gamma)$ and $\lambda_{*}(\alpha,\gamma)$ are well defined for any admissible pair $(\alpha,\gamma)$.
To cut down  notations we won't always indicate $\alpha$ and $\gamma$. \vskip0.1in

Let now give the following standard existence result.\vskip0.1in

\begin{theorem}\label{2.1}
For every $0\leq f\in L^{1}(\Omega)$ there exists a unique $0\leq
u\in L^{1}(\B)$ which satisfies
$$
\int_{\B}u\Delta^{2}\varphi dx=\int_{\B}f\varphi dx
$$
for all  $\varphi \in C^{4}(\bar{\B})\cap H_{0}^{2}(\B)$.
\end{theorem}
The proof is standard, please see \cite{Gu1}, here we omit it. From
this Theorem, we immediately have the following result. \vskip0.1in

\begin{proposition}\label{P2.1}
 Assume the existence of a weak super-solution $U$ of $(2.1)_{\lambda,\alpha,\gamma}$.
 Then there exists a weak solution $u$ of $(2.1)_{\lambda,\alpha,\gamma}$ so that $\alpha\leq u\leq U$ a.e in $\B$.\vskip0.1in
\end{proposition}
For the sake of completeness, we include a brief proof here, which be called  \textquotedblleft weak" iterative scheme:
$u_{0}=U$ and (inductively) let $u_{n}, n\geq 1$, be the solution of
$$
\int_{\B}(u_{n}-\Phi)\Delta^{2}\varphi=\lambda \int_{\B}\frac{\varphi}{(1-u_{n-1})^{p}}\ \
\forall \varphi \in C^{4}(\bar{\B})\cap H_{0}^{2}(\B),
$$
given by Theorem 2.1. Since $\alpha$ is a sub-solution of $(2.1)_{\lambda,\alpha,\gamma}$, inductively it is easily shown by Lemma 2.2
that $\alpha\leq u_{n+1}\leq u_{n}\leq U$ for every $n\geq 0$. Since
$$
(1-u_{n})^{-p}\leq (1-U)^{-p}\ \ \in L^{1}(\B),
$$
by Lebesgue Theorem the function $u=\lim_{n\to\infty}u_{n}$ is a weak solution of $(2.1)_{\lambda,\alpha,\gamma}$
so that $\alpha\leq u\leq U$. \vskip0.1in

In particular, for every $\lambda\in (0,\lambda_{*})$, we can find a weak solution of $(2.1)_{\lambda,\alpha,\gamma}$.
In the same range of $\lambda's$, this is still true for regular weak solutions as shown in the following lemma. \vskip0.1in

\begin{lemma}\label{L2.4}
Let $(\alpha, \gamma)$ be an admissible pair and $u$ be a weak
solution of $(2.1)_{\lambda,\alpha,\gamma}$. Then, there exists a
regular solution for every $0<\mu< \lambda$.
\end{lemma}
 {\bf Proof.} Let $\epsilon\in (0,1)$ be given and let
$\bar{u}=(1-\epsilon)u+\epsilon\Phi$, where $\Phi$ is given in
(2.2). by lemma 2.2 $\sup_{\B}\Phi< \sup_{\B} u\leq 1$. Hence
$$
\sup_{\B} \bar{u}\leq (1-\epsilon)+\epsilon \sup_{\B}\Phi<1,\ \  \inf_{\B} \bar{u}\geq (1-\epsilon)\alpha+\epsilon\inf_{\B}\Phi=\alpha
$$
\begin{eqnarray*}
\int_{\B}(\bar{u}-\Phi)\Delta^{2}\varphi&=&(1-\epsilon)\int_{\B}(u-\Phi)\Delta^{2}\varphi
=(1-\epsilon)\lambda\int_{\B}\frac{\varphi}{(1-u)^{p}}\\
&=&(1-\epsilon)^{p+1}\lambda \int_{\B}\frac{\varphi}{(1-\bar{u}+\epsilon(\Phi-1))^{p}}\geq
(1-\epsilon)^{p+1}\lambda\int_{\B}\frac{\varphi}{(1-\bar{u})^{p}}.
\end{eqnarray*}
Note that $0\leq (1-\epsilon)(1-u)=1-\bar{u}+\epsilon(\Phi-1)<1-\bar{u}$. So $\bar{u}$ is a weak super-solution of
$(2.1)_{(1-\epsilon)^{p+1}\lambda,\alpha,\gamma}$ such that $\sup_{\B}\bar{u}<1$. By Lemma 2.2 we get the existence of a weak solution
of $(2.1)_{(1-\epsilon)^{p+1}\lambda,\alpha,\gamma}$ so that $\alpha\leq \omega\leq\bar{u}$. In particular, $\sup_{\B}\bar{u}<1$ and $\omega$
is a regular weak solution. Since $\epsilon\in (0,1)$ is arbitrarily chosen, the proof is done.\vskip0.1in

Now we recall some basic facts about the minimal branch\vskip0.1in

\begin{theorem} \label{T2.2}
$\lambda^{*}\in (0,+\infty)$ and the following holds:\vskip0.1in

 1. For each $0<\lambda<\lambda^{*}$ there exists a regular and minimal solution $u_{\lambda}$ of $(2.1)_{\lambda,\alpha,\gamma}$;

 2. For each $x\in \B$ the map $\lambda\rightarrow u_{\lambda}(x)$ is strictly increasing on $(0,\lambda^{*})$;

 3. For $\lambda>\lambda^{*}$ there are no weak solutions of $(2.1)_{\lambda,\alpha,\gamma}$.\vskip0.1in
\end{theorem}
The proof is standard, see \cite{Co1}, here we omit it.

\setcounter{equation}{0}
 \setcounter{section}{2}
\section{Stability of the minimal solutions}
In this section we shall show that the extremal solution is regular
in small dimensions. Let us begin with the following priori
estimates along the minimal branch $u_{\lambda}$. In order to
achieve this,  we shall need yet another notion of
\emph{$H^{2}$}($\B$)- weak solutions, which is an intermediate class
between classical and weak solutions. \vskip0.1in

\begin{definition}\label{D3.1}
We say that $u$ is a \emph{$H^{2}$}($\B$)- weak solution of
$(2.1)_{\lambda,\alpha,\beta}$ if $u-\Phi \in $
\emph{$H_{0}^{2}$}($\B$), $\alpha\leq u\leq 1 \in \B$,
$\frac{1}{(1-u)^{p}}\in L^{1}(\B)$ and if
$$
\int_{\B}\Delta u\Delta\phi=\lambda\int_{\B}\frac{\phi}{(1-u)^{p}},\ \ \ \forall \phi\in C^{4}(\bar{\B})\cap H_{0}^{2}(\B),
$$
where $\Phi$ is given in (2.2). We say that $u$ is a \emph{$H^{2}$}($\B$)- weak super-solution (resp. \emph{$H^{2}$}($\B$)-
weak sub-solution) of $(2.1)_{\lambda,\alpha,\beta}$ if for $\phi\geq0$ the equality is replaced with $\geq$ (resp.$\leq$)
 and $u\geq \alpha$ (resp. $\leq$), $\frac{\partial u}{\partial n}\leq \beta$ (resp. $\geq$) on $\partial \B$.\vskip0.1in
\end{definition}

\begin{theorem}\label{T3.1}
Suppose that $(\alpha, \gamma)$ is an admissible pair.

\vskip0.1in\noindent
 1. The minimal solution $u_{\lambda}$ is stable, and is the unique
semi-stable $H^{2}(\B)$ weak solution of
$(2.1)_{\lambda,\alpha,\gamma}$;

\vskip0.1in\noindent
 2. The function $u^{*}:=\lim_{\lambda\to \lambda^{*}}u_{\lambda}$
is a well-defined semi-stable  $H^{2}(\B)$ weak solution of
$(2.1)_{\lambda^{*},\alpha,\gamma}$;

\vskip0.1in\noindent
 3.  $u^{*}$ is the unique $H^{2}(\B)$ weak solution of
$(2.1)_{\lambda^{*},\alpha,\gamma}$, and when $u^{*}$ is classical
solution, then $\mu_{1}(u^{*})=0$;

\vskip0.1in\noindent
 4.  If $v$ is a singular, semi-stable $H^{2}(\B)$ weak solution of
$(2.1)_{\lambda,\alpha,\gamma}$, then $v=u^{*}$ and
$\lambda=\lambda^{*}$.
\end{theorem}
The main tool which we use to prove the Theorem \ref{T3.1} is the
following comparison lemma which is valid exactly in the class
$H^{2}(\B)$.\vskip0.1in

\begin{lemma}\label{L3.1}
 Let $(\alpha, \gamma)$ is an admissible pair and $u$ be a semi-stable $H^{2}(\B)$ weak solution of $(2.1)_{\lambda,\alpha,\gamma}$.
 Assume $U$ is a $H^{2}(\B)$ weak super-solution of $(2.1)_{\lambda,\alpha,\gamma}$.
 Then
 \vskip0.1in\noindent
1. $u\leq U$ a.e. in $\B$; \vskip0.1in\noindent
 2. If $u$ is a classical solution and $\mu_{1}(u)=0$ then $U=u$.
\end{lemma}
A more general version of Lemma \ref{L3.1} is available in the
following.\vskip 0.1in

\begin{lemma}\label{L3.2}
Let $(\alpha, \gamma)$ is an admissible pair and $\gamma'\leq 0$.
Let $u$ be a semi-stable  $H^{2}(\B)$ weak sub-
 solution of $(2.1)_{\lambda,\alpha,\gamma}$ with $u=\alpha'\leq \alpha, \frac{\partial u}{\partial n}=\gamma'\geq \gamma$ on $\partial\B$.
 Assume that $U$ is a $H^{2}(\B)$ weak  super-solution of $(2.1)_{\lambda,\alpha,\gamma}$ with $U=\alpha, \frac{\partial U}{\partial n}=\gamma $ on
 $\partial \B$. Then $U\geq u$ a.e. in $\B$.
 \end{lemma}
The proof of Lemma \ref{L3.1} and Lemma \ref{L3.2} are same as
\cite{Co1, Mo}, we omit it here. Also, we need some a priori
estimates along the minimal branch $u_{\lambda}$.

\begin{lemma}\label{L3.3}
Let ($\alpha,\gamma$) be an
admissible pair. Then for every $\lambda\in(0,\lambda^{*} )$, we
have
$$p\int_{\B}\frac{(u_{\lambda}-\Phi)^{2}}{(1-u_{\lambda})^{p+1}}\leq \int_{\B}\frac{u_{\lambda}-\Phi}{(1-u_{\lambda})^{p}},$$
where $\Phi$ is given by (2.2). In particular, there is a constant $C$ independent of $\lambda$ so that
$$
\int_{\B}|\Delta u_{\lambda}|^{2}dx+\int_{\B}\frac{1}{(1-u_{\lambda})^{p+1}}\leq C. \eqno(3.1)
$$
\end{lemma}
{\bf Proof.} Testing $(2.1)_{\lambda,\alpha,\gamma}$ on $u_{\lambda}-\Phi\in W^{4,2}(\B)\cap H_{0}^{2}(\B)$. We see that
$$
\lambda\int_{\B}\frac{u_{\lambda}-\Phi}{(1-u_{\lambda})^{p}}=\int_{\B}
(\Delta(u_{\lambda}-\Phi))^{2}dx\geq p \lambda\int_{\B}\frac{(u_{\lambda}-\Phi)^{2}}{(1-u_{\lambda})^{p+1}}dx
$$
in the view of $\Delta^{2}\Phi=0$. In particular, for $\delta>0$ small we have that
\begin{eqnarray*}
\int_{|u_{\lambda}-\Phi|\geq \delta}\frac{1}{(1-u_{\lambda})^{p+1}}\leq \frac{1}{\delta^{2}} \int_{|u_{\lambda}-\Phi|\geq \delta}
\frac{(u_{\lambda}-\Phi)^{2}}{(1-u_{\lambda})^{p+1}}&\leq& \frac{1}{\delta^{2}} \int_{\B}\frac{1}{(1-u_{\lambda})^{p}}\\
\leq \delta^{p-1}\int_{\B}\frac{1}{(1-u_{\lambda})^{p+1}}+C_{\delta}
\end{eqnarray*}
by means of Young's inequality. Since for $\delta$ small
$$\int_{|u_{\lambda}-\Phi|\leq \delta}\frac{1}{(1-u_{\lambda})^{p+1}}\leq C$$
for some $C>0$, we get that
$$\int_{\B}\frac{1}{(1-u_{\lambda})^{p+1}}\leq C$$
for some $C>0$ and for every $\lambda\in (0,\lambda^{*})$. By Young'
s and H$\ddot{o}$lder's inequalities, we have
\begin{eqnarray*}
\int_{\B}|\Delta u_{\lambda}|^{2}dx&=&\int_{\B}
 \Delta u_{\lambda}\Delta \Phi dx+\lambda\int_{\B}\frac{u_{\lambda}-\Phi}{(1-u_{\lambda})^{p}}dx\\
&\leq& \delta\int_{\B}|\Delta u_{\lambda}|^{2}dx+C_{\delta}
+C\big(\int_{\B}\frac{dx}{(1-u_{\lambda})^{p+1}}\big)^{\frac{p}{p+1}}.
\end{eqnarray*}
So we have
$$\int_{\B}|\Delta u_{\lambda}|^{2}dx+\int_{\B}\frac{dx}{(1-u_{\lambda})^{p+1}}\leq C$$
where $C$ is absolute constant. \vskip 0.2in

{\bf Proof of the Theorem \ref{T3.1}} \ \ (1) \  \ Since
$\|u_{\lambda}\|_{\infty}<1$, the infimum  defining
$\mu_{1}(u_{\lambda})$ is achieved at a first eigenfunction for
every $\lambda \in (0,\lambda^{*})$. Since $\lambda\mapsto
u_{\lambda}(x)$ is increasing for every $x\in \B$, it is easily seen
that $\lambda\rightarrow \mu_{1}(u_{\lambda})$ is a decreasing and
continuous function on $(0,\lambda^{*})$. Define
$$
\lambda_{**}:=\sup \{0<\lambda<\lambda^{*}: \mu_{1}(u_{\lambda})>0\}.
$$
We have that $\lambda_{**}=\lambda^{*}$. Indeed, otherwise we would
have $\mu_{1}(u_{\lambda_{**}})=0$, and for every $\mu\in
(\lambda_{**}, \lambda^{*}), u_{\mu}$ would be a classical
super-solution of $(2.1)_{\lambda_{**},\alpha,\gamma}$. A
contradiction arises since Lemma \ref{L3.1} implies
$u_{\mu}=u_{\lambda_{**}}$. Finally, Lemma \ref{L3.1} guarantees the
uniqueness in the class of semi-stable $H^{2}(\B)$ weak
solutions.\vskip 0.1in

(2)\  \  It follows from (3.1) that $u_{\lambda}\to u^{*}$ in a pointwise sense and weakly in $H^{2}(\B)$,
and $\frac{1}{1-u^{*}} \in L^{p+1}$.
In particular, $u^{*}$ is a $H^{2}(\B)$ weak solution of $(2.1)_{\lambda_{**},\alpha,\gamma}$
which is also semi-stable as the limiting function of
the semi-stable solutions $\{u_{\lambda}\}$.\vskip 0.1in

(3)\  \  Whenever $\|u^{*}\|_{\infty}< 1$, the function $u^{*}$ is a
classical solution, and by the Implicit Function Theorem we have
that $\mu_{1}(u^{*})=0$ to prevent the continuation of the minimal
branch beyond $\lambda^{*}$. By Lemma \ref{L3.1}, $u^{*}$ is then
the unique $H^{2}(\B)$ weak solution of
$(2.1)_{\lambda_{*},\alpha,\gamma}$.\vskip 0.1in

(4)\ \  If $\lambda<\lambda^{*}$, we get by uniqueness that $v=u_{\lambda}$. So $v$ is not singular and a contradiction arises.
 Since $v$ is a semi-stable $H^{2}(\B)$ weak solution of $(2.1)_{\lambda_{*},\alpha,\gamma}$ and
$u^{*}$ is a $H^{2}(\B)$ weak super-solution of
$(2.1)_{\lambda_{*},\alpha,\gamma}$, we can apply Lemma \ref{L3.1}
to get $v\leq u^{*}$ a.e. in $\Omega$. Since $u^{*}$ is also a
semi-stable solution, we can reverse the roles of $v$ and $u^{*}$ in
Lemma \ref{L3.1} to see that $v\geq u^{*}$ a.e. in $\B$. So equality
$v=u^{*}$ holds and the proof is done

\setcounter{equation}{0}
 \setcounter{section}{3}
\section{Regularity of the extremal solution  and the Proof of the Theorem 1.1 (i)}

In this section we first show that the extremal solution is regular
in small dimensions by the uniformly bounded of $u_{\lambda}$ in
$H^{2}_{0}(\B)$.  Now we give the proof of Theorem \ref{T1.1} (i).
\vskip 0.1in

{\bf Proof of Theorem \ref{T1.1} (i).} As already observed, estimate
(3.1) implies that $f(u^{*})=(1-u^{*})^{-p}\in
L^{\frac{p+1}{p}}(\B)$. Since $u^{*}$ is radial and radially
decreasing. We need to show that $u^{*}(0)<1$ to get the regularity
of $u^{*}$. In fact, on the contrary suppose that $u^{*}(0)=1$. By
the standard elliptic regularity theory shows that $u^{*}\in
W^{4,\frac{p+1}{p}}$. By the Soblev imbedding theorem, i.e.
$W^{4,\frac{p+1}{p}}\hookrightarrow C^{m} (0<m\leq 4-\frac{pn}{p+1},
1\leq n\leq 4)$. We have $u^{*}$ is a Lipschitz function in $\B$ for
$1\leq n\leq 3$.

Now suppose $u^{*}(0)=1$ and $1\leq n \leq 2$. Since
$$
\frac{1}{1-u^{*}}\geq\frac{C}{|x|} \ \ \ \mbox{in} \ \ \B
$$
for some $C>0$. One see that
$$
+\infty=C^{p+1}\int_{\B}\frac{1}{|x|^{p+1}}\leq \int_{\B}\frac{1}{(1-u^{*})^{p+1}}<+\infty.
$$
A contradiction arises and hence $u^{*}$ is regular for $1\leq
n\leq2$.

For $n=3$, by the Sobolev imbedding theorem, we have $ u^{*}\in
C^{\frac{p+4}{p+1}}(\bar{B})$, if $\frac{p+4}{p+1}\geq2$, then $u^{*}(0)=1,
Du^{*}(0)=0$ and
$$
|Du^{*}(\varepsilon)-Du^{*}(0)|\leq M|\varepsilon|\leq M|x|
$$
where $0<|\varepsilon|<|x|$. Thus
$$
|u(x)-u(0)|\leq |Du(\varepsilon)||x|\leq M|x|^{2}.
$$
This inequality shows that
$$
+\infty=C^{p+1}\int_{\B}\frac{1}{|x|^{2(p+1)}}\leq \int_{\B}\frac{1}{(1-u^{*})^{p+1}}<+\infty.
$$
A contradiction arises and hence $u^{*}$ is regular for $n=3$; if $\frac{p+4}{p+1}<2$, then
$$
|Du(\varepsilon)-Du(0)|\leq M|\varepsilon|^{\frac{4}{p-1}-1}\leq M|x|^{\frac{3}{p+1}}
$$
where $0<|\varepsilon|<|x|$. Thus
$$
|u(x)-u(0)|\leq |Du(\varepsilon)||x|\leq M|x|^{\frac{4+p}{p+1}},
$$
and a contradiction  is obtained as above.

For $n=4$, by the Sobolev imbedding theorem, we have $ u^{*}\in
C^{\frac{4}{p+1}}(\bar{\B})$.  If $1<\frac{4}{p+1}<2$, then $u^{*}(0)=1,
Du^{*}(0)=0$ and
$$
|Du(\varepsilon)-Du(0)|\leq M|\varepsilon|^{\frac{4}{p+1}-1}\leq M|x|^{\frac{4}{p+1}-1}
$$
where $0<|\varepsilon|<|x|$. Thus
$$
|u(x)-u(0)|\leq |Du(\varepsilon)||x|\leq M|x|^{\frac{4}{p+1}}.
$$
If $\frac{4}{p+1}\leq 1$, then $u^{*}$ is a H$\ddot{o}$lder's
continues and
$$
1-u^{*}(x)\leq M|x|^{\frac{4}{p+1}},
$$
and we obtain a contradiction as above.
\vskip 0.1in

Now we give the point estimate of singular extremal solution for dimensions $n\geq5$.\vskip 0.1in

\begin{theorem}\label{T4.1}
Let $n\geq 5$ and $(u^{*}, \lambda^{*})$ be the extremal pair of
$(1.1)_{\lambda}$, when $u^{*}$ is singular, then
$$
1-u^{*}\leq C_{0} |x|^{\frac{4}{p+1}},
$$
where $ C_{0}:=(\lambda^{*}/K_{0})^{\frac{1}{p+1}}$ and
$K_{0}:=\frac{8(p-1)}{(p+1)^{2}}\left[n-\frac{2(p-1)}{p+1}\right]\left[n-\frac{4p}{p+1}\right]$.
\end{theorem}
In order to prove the Theorem \ref{T4.1}, we need the lower bounds
of $\lambda^{*}$ and state as follows.\vskip0.2in

\begin{lemma}\label{L4.1}
$\lambda^{*}$ satisfies the following lower bounds for $n\geq 4$
$$ \lambda^{*}\geq K_{0}$$
where
$K_{0}:=\frac{8(p-1)}{(p+1)^{2}}\left[n-\frac{2(p-1)}{p+1}\right]\left[n-\frac{4p}{p+1}\right]$
.\end{lemma}
{\bf Proof.} the proof is standard, here we include the
proof for the sake of completeness. Notice that for $n\geq 4$ the
function $\bar{u}=1-|x|^{\frac{4}{p+1}}$ satisfies
$$\frac{1}{(1-\bar{u})^{p}}\in L^{1}(\B)$$ and $\bar{u}$ is a weak
solution of
$$
\Delta^{2}\bar{u}= \frac{K_{0}}{(1-\bar{u})^{p}},
$$
and $\bar{u}(1)=0=u_{\lambda}(1); \frac{\partial u_{\lambda}}{\partial n}(1)\geq  \frac{\partial\bar{u}_{\lambda}}{\partial n} (1)$.
Therefore, $\bar{u}$ turns out to be a weak super-solution of
$(1.1)_{\lambda}$ provided $\lambda\leq K_{0}$. Thus
necessarily, we have
$$
\lambda^{*}=\lambda_{*}\geq  K_{0}.
$$
The proof is done.\vskip 0.2in

{\bf Proof of Theorem \ref{T4.1}.}\ \ First note that Lemma 4.1
gives the lower bound:
$$
\lambda^{*}\geq K_{0}.
$$
For $\delta>0$, we define $u_{\delta}(x):=1-C_{\delta}|x|^{\frac{4}{p+1}}$ with
$C_{\delta}:=(\frac{\lambda^{*}}{K_{0}}+\delta)^{\frac{1}{p+1}}>1$.  Since $n\geq5$, we have that
$u_{\delta}\in H_{loc}^{2}(\R^{n}), \frac{1}{1-u_{\delta}}\in L_{loc}^{3}(\R^{n})$ and $u_{\delta}$ is a $H^{2}$-weak solution of
$$
\Delta^{2}u_{\delta}=\frac{\lambda^{*}+\delta K_{0}}{(1-u_{\delta})^{p}}\ \ \mbox{in} \ \ \R^{n}
$$
We claim that $u_{\delta}\leq u^{*}$ in $\B$, which will finish the proof by just letting $\delta\to 0$.

Assume by contradiction that the set
$$
\Gamma:=\{r\in (0,1): u_{\delta}(r)>u^{*}(r)\}
$$
is non-empty, and let $r_{1}=\sup\Gamma.$ Since
$$
u_{\delta}(1)=1-C_{\delta}<0=u^{*}(1),
$$
we have that $0<r_{1}<1$ and one infers that
$$
\alpha:=u^{*}(r_{1})=u_{\delta}(r_{1}),\ \  \beta=(u^{*})'(r_{1})\geq u_{\delta}'(r_{1}).
$$
Setting $u_{\delta,r_{1}}(r)=r_{1}^{-\frac{4}{p+1}}(u_{\delta}(r_{1}r)-1)+1$, we easily see that the function
$u_{\delta,r_{1}}(r)$ is a $H^{2}(\B)$- weak super-solution of $(2.1)_{\lambda^{*}+\delta K_{0},\alpha',\beta'}$, where
$$
\alpha':=r_{1}^{-\frac{4}{p+1}}(u^{*}(r_{1}r)-1)+1, \ \ \beta':=r_{1}^{\frac{p-3}{p+1}}\beta.
$$
Similarly, define $u_{r_{1}}^{*}=r_{1}^{-\frac{4}{p+1}}(u^{*}(r_{1}r)-1)+1$, we have $u_{r_{1}}^{*}$ is singular semi-stable
 $H^{2}(\B)$-weak solution of $(2.1)_{\lambda^{*},\alpha',\beta'}$.

Now we claim that $(\alpha',\beta')$ is an admissible pair. Since $u^{*}$ is radially decreasing, we have that
$\beta'\leq0$. Define the function
$$
\omega(r):=(\alpha'-\frac{\beta'}{2})+\frac{\beta'}{2}|x|^{2}+\gamma(x),$$
where $\gamma(x)$ is a solution of $\Delta^{2} \gamma=\lambda^{*}$ in $\B$ with $\gamma=\partial_{v} \gamma=0$ on $\partial \B$.
 Then $\omega$ is a classical solution of
$$
\left\{ \begin{array}{lllllll}
\Delta^{2}\omega=\lambda^{*} & \mbox{in} \ \  \B \\
\omega=\alpha', \partial_{v}\omega=\beta' & \mbox{on} \ \  \partial \B.
\end{array}
\right.
$$
Since $\frac{\lambda^{*}}{(1-u^{*}_{r_{1}})^{p}}\geq \lambda^{*}$, by Lemma 2.1 we have
$$
u^{*}_{r_{1}}\geq \omega \ \ \mbox{a.e.} \ \ \mbox{in}\ \  \B
$$
Since $\omega(0)=\alpha'-\frac{\beta'}{2}+\gamma(0)$ and $\gamma(0)>0$, we have
$$
\alpha'-\frac{\beta'}{2}<1
$$
So $(\alpha',\beta')$ is an admissible pair and by Theorem 3.1 (4) we get that $(u^{*}_{r_{1}}, \lambda^{*})$
coincides with the extremal pair of $(2.1)_{\lambda,\alpha',\beta'}$ in $\B$.

Since $(\alpha',\beta')$ is an admissible pair and
$u_{\delta,r_{1}}$ is a $H^{2}(\B)$-weak super-solution of
$(2.1)_{\lambda^{*}+\delta K_{0},\alpha',\beta'}$. We get from
Proposition \ref{P2.1}, the existence of a weak solution of
$(2.1)_{\lambda^{*}+\delta K_{0},\alpha',\beta'}$. Since
$$
\lambda^{*}+\delta K_{0}>\lambda^{*},
$$
we contradict the fact that $\lambda^{*}$ is the extremal parameter of $(2.1)_{\lambda,\alpha',\beta'}$.\vskip0.1in

Thanks to the lower estimate on $u^{*}$, we get the following  result.\vskip0.1in

\begin{corollary}\label{C4.1}
 In any dimension $n\geq 1$, we have
$$
\lambda^{*}>K_{0}= \frac{8(p-1)}{(p+1)^{2}}\left[n-\frac{2(p-1)}{p+1}\right]\left[n-\frac{4p}{p+1}\right].
$$
\end{corollary}
\vskip0.1in {\bf Proof.}  The function
$\bar{u}:=1-|x|^{\frac{4}{p+1}}$ is a $H^{2}(\B)$-weak solution of
$(2.1)_{K_{0},0,-\frac{4}{p+1}}$. If by contradiction
$\lambda^{*}=K_{0}$, then $\bar{u}$ is a $H^{2}(\B)$-weak
super-solution of $(1.1)_{\lambda}$ for every $\lambda\in
(0,\lambda^{*})$. By Lemma \ref{L3.1} we get
 that $u_{\lambda}\leq\bar{u}$ for all $\lambda<\lambda^{*}$, and then $u^{*}\leq \bar{u}$ a.e. in $\B$.

If $1\leq n\leq 4$, $u^{*}$ is then regular by Theorem (i). By
Theorem \ref{T3.1} (3) there holds $\mu_{1}(u^{*})=0$. By Lemma 3.2
then yields that $u^{*}=\bar{u}$, which is a contradiction since
then $u^{*}$ will not satisfy the boundary conditions.

If now $n\geq 5$ and $\lambda^{*}=K_{0}$, then $C_{0}=1$ in Theorem
\ref{T4.1}, and we then have $u^{*}\geq\bar{u}$. It means again that
$u^{*}=\bar{u}$, a contradiction that completes the
proof.\vskip0.1in

In what follows, we will show that the extremal solution $u^{*}$ of $(1.1)_{\lambda}$ in dimensions $n\geq 13$ is singular.\vskip 0.2in

\setcounter{equation}{0}
 \setcounter{section}{4}
\section{The extremal solution is singular for $n\geq 13$}
We prove in this section that the extremal solution is singular for
$n\geq 13$ and $p$ large enough.
 For that we will need a a suitable Hardy-Rellich type inequality which was established by Ghoussoub-Moradifam in \cite{Gh}.
  As in the previous section $(u^{*},\lambda^{*})$ denotes the extremal pair of $(2.1)_{\lambda^{*},0,0}$\vskip 0.1in

\begin{lemma}\label{L5.1}
 Let $n\geq 5$ and $\B$ be the unit ball in $\R^{n}$. Then there exists $C>0$, such that
 the following improved Hardy-Rellich inequality holds for all $\varphi\in H_{0}^{2}(\B)$:
$$
\int_{\B}(\Delta \varphi)^{2}dx\geq \frac{n^{2}(n-4)^{2}}{16}\int_{\B}\frac{\varphi^{2}}{|x|^{4}}dx+C\int_{\B}\varphi^{2}dx
$$
\end{lemma}

\begin{lemma}\label{L5.2}
 Let $n\geq 5$ and $\B$ be the unit ball in $\R^{n}$. Then the following improved Hardy-Rellich
inequality holds for all $\varphi\in H_{0}^{2}(\B)$:
\begin{eqnarray*}
\int_{\B}(\Delta \varphi)^{2} dx &\geq& \frac{(n-2)^{2}(n-4)^{2}}{16}\int_{\B}\frac{\varphi^{2}dx}{(|x|^{2}-0.9|x|^{\frac{n}{2}+1})(|x|^{2}-|x|^{\frac{n}{2}})}\\
&+& \frac{(n-1)(n-4)^{2}}{4}\int_{\B}\frac{\varphi^{2}dx}{|x|^{2}(|x|^{2}-|x|^{\frac{n}{2}})}. \hspace*{4.7cm} (5.0)
\end{eqnarray*}
As a consequence, the following improvement of the classical Hardy-Rellich inequality holds:
$$
\int_{\B}(\Delta \varphi)^{2} dx\geq \frac{n^{2}(n-4)^{2}}{16}\int_{\B}\frac{\varphi^{2}}{|x|^{2}(|x|^{2}-|x|^{\frac{n}{2}})}.
$$
\end{lemma}

\begin{lemma}\label{L5.3}
  If $n\geq 13$, then $u^{*}\leq 1-|x|^{\frac{4}{p+1}}$.
\end{lemma}
{\bf Proof.}\  Recall from Corollary 4.1 that $K_{0}<\lambda^{*}$. Let $\bar{u}=1-|x|^{\frac{4}{p+1}}$, we now claim that
$u_{\lambda}\leq \bar{u}$ for all $\lambda\in (K_{0},\lambda^{*})$. Indeed, fix such a $\lambda$ and assume by contradiction that
$$
R_{1}:=\inf\{0\leq R\leq 1: u_{\lambda}<\bar{u}\ \  \mbox{in the interval } (R,1)\}>0.
$$
From the boundary condition, one has that $u_{\lambda}<\bar{u}(r)$ as $r\to 1^{-}$. Hence, $0<R_{1}<1,
\alpha:=u_{\lambda}(R_{1})=\bar{u}(R_{1})$ and $\beta:=u'_{\lambda}(R_{1})\leq\bar{u}'(R_{1}).$
The same as the proof of Theorem 4.1, we have $u_{\lambda}\geq \bar{u}$ in $\B_{R_{1}}$ and a contradiction arises in view of the
 fact that $\lim_{x\to 0}\bar{u}(x)=1$ and $\|u_{\lambda}\|_{\infty}<1$. It follows that
 $u_{\lambda}\leq \bar{u}$ in $\B$ for every $\lambda \in (K_{0}, \lambda^{*})$ and in particular
 $u^{*}\leq \bar{u}$ in $\B$.

\vskip 0.1in
\begin{lemma}\label{L5.4}
 Let $n\geq 13$.  Suppose there
exists $\lambda'>0$ and a singular radial function $\omega(r) \in
H^{2}(\B)$ with $\frac{1}{1-\omega}\in
L_{loc}^{\infty}(\bar{\B}\setminus\{0\})$ such that
$$\left\{
\begin{array}{lllllll}
\Delta^{2}\omega\leq\frac{\lambda'}{(1-\omega)^{p}}\ \ \ \ \ \ \ \ \ \  \mbox{for}\  \ 0<r<1,\\
\omega(1)=0, \ \ \ \omega'(1)=0,\\
\end{array}
\right.
\eqno(5.1)
$$
and
$$
p\beta\int_{\B}\frac{\varphi^{2}}{(1-\omega)^{p+1}}\leq\int_{\B}(\Delta\varphi)^{2}\ \ \mbox{ for all}\  \varphi\in  H_{0}^{2}(\B)
\eqno(5.2)
$$
\ \ \ \ \ 1. If $\beta\geq \lambda'$, then $\lambda^{*}\leq \lambda'$.

2. If either $\beta>\lambda'$ or $\beta=\lambda'=\frac{H_{n}}{p}$,
then the extremal solution $u^{*}$ is necessarily singular.
\end{lemma}
{\bf Proof.} (1).  First, note that (5.2) and $\frac{1}{1-\omega}\in
L_{loc}^{\infty}(\bar{\B}\setminus\{0\})$ yield to
$$
\frac{1}{1-\omega}\in L^{1}(\B).
$$
At the same time, (5.1) implies that $\omega(r)$ is a $H^{2}(\B)$-
weak stable sub-solution of $(1.1)_{\lambda'}$. If now
$\lambda'<\lambda^{*}$, then by Lemma \ref{L3.2}, we have
$$
\omega(r)<u_{\lambda'},
$$
which ia a contradiction since $\omega$ is singular while $u_{\lambda'}$ is regular.\vskip 0.1in

(2) Suppose first that $\beta=\lambda'=\frac{H_{n}}{p}$ and that $n\geq 13$. Since by part (1) we have
 $\lambda^{*}\leq \frac{H_{n}}{p}$, we get from Lemma \ref{L5.3} and improved Hardy-Rellich inequality  that there exists $C>0$ so that
 for all $\phi\in H_{0}^{2}(\B)$
 $$
 \int_{\B}(\Delta \phi)^{2}-p\lambda^{*}\int_{\B}\frac{\phi^{2}}{(1-u^{*})^{p+1}}\geq
  \int_{\B}(\Delta \phi)^{2}-H_{n}\int_{\B}\frac{\phi^{2}}{|x|^{4}}\geq C\int_{\B}\phi^{2}.
  $$
It follows that $\mu_{1}(u^{*})>0$ and $u^{*}$ must therefore be singular since otherwise, one could use
the Implicit Function Theorem to continue the minimal branch beyond $\lambda^{*}$.

Suppose now that $\beta>\lambda'$, and let $\frac{\lambda'}{\beta}<\gamma<1$ in such a way that
$$
\alpha:=(\frac{\gamma \lambda^{*}}{\lambda'})^{\frac{1}{p+1}}<1.
$$
Setting $\bar{\omega}:=1-\alpha(1-\omega)$, we claim that
$$
u^{*}\leq \bar{\omega} \ \ \ \mbox{in}\ \ \ \B. \eqno(5.3)
$$
Note that by the choice of $\alpha$ we have $\alpha^{p+1}\lambda'<\lambda^{*}$, and therefore to prove (5.3)
 it suffices to show that for $\alpha^{p+1}\lambda'\leq \lambda<\lambda^{*}$, we have $u_{\lambda}\leq \bar{\omega}$ in $\B$.
 Indeed, fix such $\lambda$ and note that
$$
\Delta^{2}\bar{\omega}=\alpha \Delta^{2}\omega\leq\frac{\alpha \lambda'}{(1-\omega)^{p}}
=\frac{\alpha^{p+1}\lambda'}{(1-\bar{\omega})^{p}}\leq \frac{\lambda}{(1-\bar{\omega})^{p}}.
$$
Assume that $u_{\lambda}\leq \bar{\omega}$ dose not hold in $\B$, and consider
$$
R_{1}:=\sup\{0\leq R\leq 1| u_{\lambda}(R)>\bar{\omega}(R)\}>0
$$
Since $\bar{\omega}(1)=1-\alpha>0=u_{\lambda}(1)$, we then have
$$
R_{1}<1, u_{\lambda}(R_{1})=\bar{\omega}(R_{1})\  \mbox{and}\  u_{\lambda}'(R_{1})\leq\bar{\omega}'(R_{1}).
$$
Introduce, as in the proof of the Theorem 4.1, the functions $u_{\lambda, R_{1}}$ and $\bar{\omega}_{R_{1}}$.
 We have that $u_{\lambda,  R_{1}}$ is a classical solution of $(2.1)_{\lambda, \alpha',\beta'}$, where
$$
\alpha':=R_{1}^{-\frac{4}{p+1}}(u_{\lambda}(R_{1})-1)+1, \beta':=R_{1}^{\frac{p-3}{p+1}}u_{\lambda}'(R_{1}).
$$
Since $\lambda<\lambda^{*}$ and then
$$
\frac{p\lambda}{(1-\bar{\omega})^{p+1}}\leq\frac{p\lambda^{*}}{\alpha^{p+1}(1-\omega)^{p+1}}
<\frac{p\beta}{(1-\omega)^{p+1}},
$$
by (5.2) $\bar{\omega}_{R_{1}}$ is a stable $H^{2}(\B)$-weak
sub-solution of $(2.1)_{\lambda,\alpha',\beta'}$. By Lemma
\ref{L3.2}, we deduce that $u_{\lambda}\geq \bar{\omega}$ in
$\B_{R_{1}}$ which is impossible, since $\bar{\omega}$ is singular
while $u_{\lambda}$ is regular. This establishes claim (5.3) which,
combined with the above inequality, yields
$$
\frac{p\lambda^{*}}{(1-u^{*})^{p+1}}\leq\frac{p\lambda^{*}}{\alpha^{p+1}(1-\omega)^{p+1}}
<\frac{p\beta}{(1-\omega)^{p}},
$$
and thus
$$
\inf_{\varphi\in C_{0}^{\infty}(\B)}\frac{\int_{\B}[(\Delta \varphi)^{2}-\frac{p\lambda^{*}\varphi^{2}}{(1-u^{*})^{p+1}}]dx}{\int_{\B}\varphi^{2}dx}>0.
$$
This is not possible if $u^{*}$ is a smooth function, since otherwise, one could use the Implicit function Theorem
 to continue the minimal branch beyond $\lambda^{*}$.

\vskip0.1in

{ \bf Proof Theorem \ref{T1.1} (ii).}

 Now we prove that $u^{*}$ is a singular solution of $(1.1)_{\lambda^{*}}$ for $n\geq 13$, in order to achieve this,
 we shall find  a singular $H^{2}(\B)$ weak sub-solution of $(1.1)_{\lambda'}$, denote by $\omega_{m}(r)$,
 which is stable, according to the Lemma \ref{L5.4}.\vskip 0.1in

Choosing
$$\omega_{m}(r)=1-a_{1}r^{\frac{4}{p+1}}+a_{2}r^{m},\quad K_{0}= \frac{8(p-1)}{(p+1)^{2}}\left[n-\frac{2(p-1)}{p+1}\right]
\left[n-\frac{4p}{p+1}\right],$$
since $\omega(1)=\omega'(1)=0$, we have
$$
a_{1}=\frac{m}{m-\frac{4}{p+1}};\quad
a_{2}=\frac{\frac{4}{p+1}}{m-\frac{4}{p+1}}.
$$
For any $m$ fixed, when $p\to \infty$, we have
$$a_1 = 1 + \frac{4}{(p+1)m} + o(p^{-1}) \quad \mbox{and}\quad \mbox a_2 = a_1 - 1 = \frac{4}{(p+1)m} + o(p^{-1})$$
and
$$K_{0} = \frac{8(n-2)(n-4)}{p} + o(p^{-1}).$$
 Note that
\begin{eqnarray*}
\frac{\lambda'K_{0}}{(1-\omega_{m}(r))^{p}}&-&\Delta^{2}\omega_{m}(r)=\frac{\lambda'K_{0}}{(1-\omega_{m}(r))^{p}}
-a_{1}K_{0}r^{-\frac{4p}{p+1}}- K_{1}r^{m-4}\\
&=&\frac{\lambda'K_{0}}{(a_{1}r^{\frac{4}{p+1}}-a_{2}r^m)^{p}}-a_{1}K_{0}r^{-\frac{4p}{p+1}}
-a_{2}K_{1}r^{m-4}\\
&=&K_{0}r^{-\frac{4p}{p+1}}\left[\frac{\lambda'}{(a_{1}-a_{2}r^{m-\frac{4}{p+1}})^{p}}-a_{1}
-a_2K_{1}K_{0}^{-1}r^{\frac{4p}{p+1}+m-4}\right]\\
&=&K_{0}r^{-\frac{4p}{p+1}}\left[\frac{\lambda'}{(a_{1}-a_{2}r^{m-\frac{4}{p+1}})^{p}}-a_{1}
-a_2 K_{1}K_{0}^{-1}r^{m -\frac{4}{p+1}}\right]\\
&=&\frac{K_{0}r^{-\frac{4p}{p+1}}}{(a_{1}-a_{2}r^{m-\frac{4}{p+1}})^p}\left[\lambda' - H(r^{m -\frac{4}{p+1}})\right] \hspace*{4.7cm} (5.4)\\
  \end{eqnarray*}
with
$$H(x) = (a_1+a_2x)^p\left[a_{1}+a_2K_{1}K_{0}^{-1}x\right], \  K_{1}=m(m-2)(m+n-2)(m+n-4) \eqno(5.5)$$

(1) Let $m=2, n\geq 32$, then we can prove that
 $$
\sup_{[0, 1]} H(x) = H(0) = a_1^{p+1}\longrightarrow e^{2}\ \  \mbox{as}\ \  p \longrightarrow +\infty. $$
So $(5.4)\geq0$ is valid as long as
$$\lambda' = e^{2}.$$
At the same time, we have (since $a_{1}-a_{2}r^{2-\frac{4}{p+1}} \geq a_1 - a_2 \geq 1$ in $[0, 1]$)
$$
\frac{n^2(n-4)^{2}}{16}\frac{1}{r^{4}} -\frac{p\beta}{r^4(a_{1}-a_{2}r^{2-\frac{4}{p+1}})^{p+1}}
 \geq r^{-4}\left[\frac{n^2(n-4)^{2}}{16} - p\beta\right].\eqno(5.6)
$$

Let $\beta=(\lambda'+\varepsilon)K_{0}$, where $\varepsilon$ is arbitrary sufficient small,  we need finally here
$$\frac{n^2(n-4)^{2}}{16} - p\beta = \frac{n^2(n-4)^{2}}{16} - p(\lambda'+\varepsilon)K_{0} >0.$$
For that, it is sufficient to have for $p \longrightarrow +\infty$
$$\frac{n^2(n-4)^{2}}{16} - 8(e^{2}+\varepsilon)(n-2)(n-4)+o(\frac{1}{p})>0.$$
So $(5.6)\geq0$ holds only for $n\geq 32$ when $p\longrightarrow +\infty$. Moreover, for $p$ large enough
$$
8e^{2}(n-2)(n-4)\int_{\B}\frac{\varphi^{2}}{(1-\omega_{2})^{p+1}}\leq
H_{n}\int_{\B}\frac{\varphi^{2}}{|x|^{4}}\leq \int_{\B}|\Delta \varphi|^{2}
$$
Thus it follows from Lemma \ref{L5.4} that $u^{*}$ is singular with
$\lambda'=e^{2}K_{0}, \beta=(e^{2}K_{0}+\varepsilon(n,p))$ and
$\lambda^{*}\leq e^{2}K_{0}$

(2) Assume $13\leq n\leq 31$. We shall show that $u=\omega_{3.5}$
satisfies the assumptions of Lemma \ref{L5.4} for each dimension
$13\leq n\leq 31$. Using Maple, for each dimension $13\leq n\leq31$
one can verify that inequality $(5.4)\geq0$ holds for the $\lambda'$
given by Table 1. Then,
 by using Maple again, we show that there exists $\beta>\lambda'$ such that
\begin{eqnarray*}
\frac{(n-2)^2(n-4)^2}{16}\frac{1}{(|x|^{2}-0.9|x|^{\frac{n}{2}+1})(|x|^{2}-|x|^{\frac{n}{2}})}\\
+\frac{(n-1)(n-4)^{2}}{4}\frac{1}{|x|^{2}(|x|^{2}-|x|^{\frac{n}{2}})}&\geq& \frac{p\beta}{(1-w_{3.5})^{p+1}}.
\end{eqnarray*}
The above inequality and and improved Hardy-Rellich inequality (5.0) guarantee that the stability condition (5.2) holds for
$\beta>\lambda'$. Hence by Lemma 5.4 the extremal solution is singular for $13\leq n\leq31$ the value of $\lambda'$ and $\beta$
are shown in Table 1.\vskip0.1in

\begin{remark}
The values of $\lambda'$ and $\beta$ in Table 1 are not optimal.
\end{remark}

$$
\begin{tabular}{|l|l|l|}
\multicolumn{3}{c}{\large Table1}\\ [5pt]
\hline
$n$ &  $\lambda'$   &$\beta$\\ \hline
31& 3.15$K_{0}$ & 4$K_{0}$ \\ \hline
30-19& 4$K_{0}$ & 10$K_{0}$ \\ \hline
18&3.19$K_{0}$&  3.22$K_{0}$ \\ \hline
17& 3.15$K_{0}$&  3.18$K_{0}$ \\ \hline
16& 3.13$K_{0}$&  3.14$K_{0}$ \\ \hline
15& 2.76$K_{0}$&  3.12$K_{0}$ \\ \hline
14&2.34$K_{0}$&  2.96$K_{0}$ \\ \hline
13&2.03$K_{0}$&  2.15$K_{0}$ \\ \hline
\end {tabular}
$$

\begin{remark}
The improved Hardy-Rellich inequality (5.0) is crucial to prove that
$u^{*}$ is singular in dimensions $n\geq 13$. Indeed by the
classical Hardy-Rellich inequality and $u:=w_{2}$, Lemma \ref{L5.4}
only implies that $u^{*}$ is singular n dimensions $n\geq 32$.
\end{remark}

\bigskip
\noindent
\textbf{Acknowledgements.} The first author is greatly indebted to
 Prof. Yi Li, his supervisor, for useful discussions. We also thank
 the anonymous referee for his valuable suggestions.
The first author is supported in part  by National Natural Science
Foundation of China (Grant No. 10971061), Natural Science Foundation
of Henan Province (Grant No. 112300410054) and Natural Science
Foundation of Education Department of Henan Province (Grant No
2011B11004). The second author is supported by the Fundamental
Research Funds for the Central Universities, Hunan University and
partially by NSFC, No 10971057.

\end{document}